\documentclass{amsart}
\usepackage{amssymb}

\makeatletter

\makeatletter
\let\cal\mathcal
\makeatother

\makeatletter\@addtoreset{equation}{section}\makeatother

\begin{document}
\bibliographystyle{plain}
\title{The Many Faces of Modern Combinatorics}
\author{Cristian Lenart}
\address{Department of Mathematics and Statistics, SUNY at Albany, Albany, NY 12222}
\email{lenart@albany.edu}
\maketitle

\begin{abstract}
This is a survey of recent developments in combinatorics. The goal is to give a big picture of (and related references for) its many interactions with other areas of mathematics, such as: group theory, representation theory, commutative algebra, geometry (including algebraic geometry), topology, probability theory, and theoretical computer science. 
\end{abstract}

Combinatorics is an area of mathematics 
that has experienced tremendous growth during the last few decades. Combinatorial mathematics was pursued since antiquity, but the subject was considered as lying at the periphery of mathematics, being mostly associated with basic counting, puzzles, and games.
As Anders Bj\"{o}rner and Richard Stanley point out in [3], the reasons for the recent growth of combinatorics come from its relationship with areas at the center stage of mathematics on the one hand, and with other disciplines (most of all computer science), on the other hand. This is a two-way relationship: combinatorics has found important applications to other areas of pure mathematics and computer
science, while certain mathematical techniques (for instance, from linear algebra, group theory, geometry, probability theory, statistical mechanics etc.) and computer experiments are now used in combinatorics. 

Let us look more closely at the factors which stimulated the recent developments in combinatorics. We already mentioned  computer science as one of these factors, due to its many mathematical problems and needs, plus potential for mathematical experiments. Combinatorics turned out to be the area of mathematics that best fit the needs of computer science. Indeed, combinatorics had an essential contribution to the formalization of structures, methods, and problems in computer science, and provided some of the tools for solving these problems. New areas, with a highly combinatorial character, were born out of this interaction, such as: combinatorial optimization, complexity theory, randomized algorithms, computational geometry, computational group theory, and computational molecular biology. On the other hand, combinatorial structures are very well suited for experiments using computer algebra systems such as {\em Mathematica} and {\em Maple}. At this moment, there are a great number of combinatorial packages for the above systems, which enable combinatorialists to gain a lot of insight into their problems, by experimenting with their structures up to reasonably high orders of magnitude. 

Beside computer science, developments within mathematics itself also influenced combinatorics. One of these developments, mentioned in [3], is that after an era where the fashion was to seek generality and abstraction, there is now much appreciation and emphasis for concrete calculations, which are the ``hard'' problems. The picture of contributions of combinatorics to mathematical calculations can be sketched by referring to three major areas of mathematics. The first two are geometry and topology, with their complicated and hard to manipulate objects. Then comes algebra, which associates various invariants and structures with these objects (such as homotopy and homology groups); they capture some of the geometric and topological information, but not all of it, of course.  Although passing to algebra is an important simplification, even calculations in this world can be very hard. This is the moment when combinatorics comes into play, since it turns out that many calculations can be done in terms of combinatorial objects. For instance, think of a summation, where the index is running from 1 to $n$; next, by contrast, think of a summation with the index running over a set of certain puzzles; it turns out that this situation arises frequently when one wants to do concrete calculations. Furthermore, if the sum consists of both positive and negative terms, one might be able to cancel the negative ones by combinatorially mapping the associated objects to the objects associated with the positive terms. Finally, one might be able to prove identities by constructing a combinatorial bijection between the objects corresponding to the terms on the two sides of the identity. Another situation is that in which one is able to find a combinatorial interpretation for positive integer cofficients appearing in several areas of mathematics. In many cases, this interpretation turns out to be very useful for computational purposes and for studying the properties of the mentioned coefficients; moreover, it sometimes points to deeper algebraic or geometric phenomena. Of course, mathematicians continue to invent more and more general algebraic and geometric structures. Unfortunately, the gap between the level of generality of the newly introduced structures and the level at which most concrete calculations can be done is wider than ever. Let us just mention  some simple geometric structures, such as Grassmannians and flag varieties, for which we still cannot do certain basic calculations (see section 7 below).

Speaking about developments within mathematics which led to the recent growth of combinatorics, we must also mention the discovery of deep connections with mainstream areas of mathematics, such as algebra, geometry, topology, and probability theory, to mention just a few. For instance, it might look like a miracle that the best frameworks for certain combinatorial problems which are easy to explain, even to the layman, are related to topology, group theory, or statistical mechanics (see some of the examples below). Indeed, tools from the mentioned areas sometimes offer more insight into the structure of combinatorial configurations. On the other hand, it turns out that certain algebraic and geometric structures are best described in combinatorial terms, and certain constructions related to them are best done with combinatorial tools. 

We have discussed the reasons for the tremendous growth of combinatorics in recent years, but we still have not explained what combinatorics is about. 
Combinatorics could be described as the study of arrangements of objects
according to specified rules. We want to know, first, whether a
particular arrangement is possible at all, and, if so, in how many
different ways it can be done. Algebraic and even probabilistic methods play an increasingly
important role in answering these questions. If we have two sets of arrangements with the same
cardinality, we might want to construct a natural bijection between them. 
We might also want to have
an algorithm for constructing a particular arrangement or all
arrangements, as well as for  computing numerical characteristics of them; in particular, we can consider optimization problems related to such arrangements. 
Finally, we might be interested in an even deeper study, by investigating the structural 
properties of the arrangements. Methods from areas such as group theory and topology are useful here, by enabling us to study symmetries of the arrangements, as well as topological properties of certain spaces associated with them, which translate into combinatorial properties. 

We are now going to have a closer look at the different faces of modern combinatorics, by briefly discussing its interaction with several areas of mathematics and computer science. We tried to emphasize the main developments in each area, including the recent ones, by presenting selected examples. However, their number is kept to a minimum and all definitions of the structures involved are at least briefly presented. 

{\bf 1. The modern face of enumerative combinatorics.} Enumeration is one of the basic and important aspects of combinatorics. In many branches of mathematics and its applications we need to know the number of ways of doing something, or at least to approximate this number. Some of the basic problems in enumerative combinatorics are to enumerate: 1) words of different types with letters from a given alphabet; 2) partitions, that is ways of breaking an object (for instance, a positive integer or a set) into smaller objects of the same kind; 3) orbits of a finite group acting on a finite set (such as colorings of the vertices, edges, or faces of a Platonic solid up to symmetry). A great number of problems of this type were solved a couple of decades or more ago. An important tool for solving them were {\em generating functions}; these are formal power series whose coefficients are the number of objects of a certain kind with a given order of magnitude, or a given numerical characteristic, also known as statistic (for instance, words of a fixed length, permutations with a given number of inversions and/or descents, colorings of an object with a fixed number of colors). It turns out that algebraic manipulations of these power series (such as multiplication, functional composition, compositional inverse etc.) correspond to combinatorial manipulations of the objects (such as breaking an object into objects of the same kind). Another remarkable fact about generating functions is that many important functions (such as the secant function, the Jacobi elliptic functions, and Bessel functions, to mention just a few) appear in this context. Closely related to generating functions are the applications of the theory of {\em orthogonal polynomials} to enumerative combinatorics. Indeed, many orthogonal polynomials occur as algebraic invariants of combinatorial structures (such as matchings polynomials of certain graphs, which are discussed in section 2). Moreover, a number of combinatorial sequences can be expressed in terms of an inner product $(\cdot,\cdot)$  on  polynomials, which corresponds to some remarkable sequence of orthogonal polynomials; one such expression for combinatorial sequences is as {\em moment sequences} $(1,x^n)$. 

Generating functions usually do not allow us to see the whole picture related to an enumeration problem. Indeed, we might be able to prove that two statistics on a finite set of objects have the same distribution; however, in order to understand the deeper reason for which this is true, we need to construct a bijection on the set of objects which maps one statistic to the other. This situation arises frequently, for instance, in the area of permutation enumeration, which is still an active area of research. 

Beside permutation enumeration, other active areas of research in enumerative combinatorics are related to objects such as {\em plane partitions} and {\em tilings}. Plane partitions are left-justified arrays of positive integers with certain properties. For instance, if the entries are weakly increasing in rows (from left to right), and strictly increasing down the columns, plus the lengths of the rows form a weakly decreasing sequence, the plane partition is called a {\em semistandard Young tableau}. The sequence of lengths of the rows of a tableau is known as its shape. See sections 4 and 7 for applications of Young tableaux. There are various nice formulas for the number of plane partitions of different types. There are also various algorithms which construct bijections between certain plane partitions and other types of combinatorial objects, such as words, nonintersecting lattice paths (such as paths with steps (1,0) and (0,1)), reduced decompositions of permutations (see section 6), and tilings.

Tilings are partitions of a geometric region into smaller regions of some specified types. For instance, think of tilings of a chessboard with dominos. We want to know first of all if a certain region is tileable; if so, we are interested in constructing a specific tiling, as well as in counting the number of tilings, or at least in approximating this number. It turns out that some of these problems can be reformulated as problems in {\em combinatorial group theory}. For instance, a necessary condition for a region to be tileable is that a certain element (determined by the boundary of the region) is equal to the identity in a certain group defined by generators and relations (both are determined by the tiles). For some regions there are nice formulas for the number of tilings. One such region is the {\em Aztec diamond}, which consists of two squares in the first row, four squares in the second row starting one square to the left of the first row, and so on up to the $n$-th row; then the diagram obtained so far is reflected about the bottom edge and added to the original one. Another example is that of certain puzzle-piece tilings of a square, which are in bijection with the generalizations of permutation matrices known as {\em alternating-sign matrices}. The number of these is also given by a simple formula. This formula was proved both combinatorially and by using methods from statistical mechanics, in particular the Yang-Baxter equation for the 6-vertex lattice model.  An alternating-sign matrix is a matrix of 0's, 1's, and $-1$'s, in which the entries in each row and column sum to 1, and the nonzero entries in each row and column alternate in sign. 

{\bf 2. Algebraic invariants and structures associated with combinatorial objects.} In the same way as algebraic structures associated with geometric and topological objects are important tools for investigating the latter, such structures play an important role in combinatorics too. We are going to mention some of them.

Let us start with polynomial invariants of graphs (a graph consists of a set of vertices and a set of edges connecting some of the vertices in pairs). There are a large number of such invariants, each enumerating certain structures on the graph. For instance, the {\em chromatic polynomial} enumerates proper colorings of the graph (that is, assignments of colors to the vertices such that no two adjacent vertices have the same color); the {\em matchings polynomial} enumerates matchings in the graph (that is, collections of edges, no two of which share a vertex). Results in this area are concerned with: 1) recurrence relations in terms of graph operations such as deletion and contraction of an edge; 2) their zeros; 3) various identities, including connections with other algebraic invariants (for instance, the chromatic polynomial of a graph can be expressed as the characteristic polynomial $-$ see the next paragraph $-$ of an associated partially ordered set); 4) interesting evaluations (for instance, the chromatic polynomial evaluated at $-1$ gives the number of orientations of the edges of the graph containing no directed cycle). Recently, there has been much interest in a symmetric function generalization of the chromatic polynomial. Some of the related results are concerned with combinatorial interpretations of the expansion of this function in various bases of the algebra of symmetric functions. 

Now we turn to invariants of partially ordered sets. A well-known one is the {\em M\"{o}bius function}; this is an integer-valued function defined on intervals as an alternating sum of the number of chains of a fixed length between the two endpoints of the corresponding interval. There are efficient methods for computing the M\"{o}bius function of certain partially ordered sets, in particular lattices with certain properties. These methods lead to nice formulas in the case of partially ordered sets of certain partitions of sets, the lattice of partitions of numbers, the lattice of subspaces of a finite-dimenional vector space over a finite field, the intersection lattice of subspace arrangements (see section 8), the Bruhat order on Coxeter groups (see section 6) etc. For certain partially ordered sets one can define  a generating function for the M\"{o}bius function, which is a polynomial called the {\em characteristic polynomial}. A nice property of this polynomial is that it factors nicely for certain lattices. Recently there was interest in a different generating function, which is a power series in infinitely many variables related to enumeration of chains in partially ordered sets. Some of the related results are concerned with sufficient conditions for this power series to be a symmetric function (in general, it satisfies a weaker condition, known as {\em quasisymmetry}). In particular, the Schur functions, introduced in section 4, can be obtained as special cases of this construction. 

Linear algebra methods also play an important role in combinatorics. For instance, one can associate with every graph on $n$ vertices its $n\times n$ {\em adjacency matrix}; an entry of this matrix is 1 or 0, depending on whether the two vertices associated with the corresponding row and column are joined by an edge or not. The eigenvalues of the adjacency matrix or of closely related matrices give information about numerical invariants such as the {\em chromatic number} (i.e., the minimum number of colors needed in a proper coloring). They also give information about certain properties of the graph, such as: strong regularity (see section 5), distance-regularity (see below), and {\em expansion} properties (measured by the ratio between the number of neighbors of vertices in a set $S$ outside $S$, and the size of $S$).

Certain combinatorial objects have algebra structures associated with them. For instance, for a large class of objects (including graphs and partially ordered sets), there is a construction of {\em Hopf algebras} (that is, algebras with a product and a dual operation called coproduct, satisfying several conditions), called {\em incidence algebras}. The coproduct describes the way in which a combinatorial object splits into smaller objects of the same kind.  The structure of incidence Hopf algebras was investigated in special cases. Incidence Hopf algebras provide the natural framework for the definition of the M\"{o}bius function, as well as for the formalization of the nineteenth century symbolic calculus (also referred to as {\em umbral calculus} or {\em operator calculus}). The main idea is that many remarkable polynomial sequences (including several sequences of orthogonal polynomials) can be defined and manipulated as elements of a natural basis for a certain Hopf algebra of polynomials in one variable, on which the dual Hopf algebra acts as an algebra of {\em shift-invariant} operators (that is, linear operators commuting with the ordinary derivative).

Finally, let us mention that certain sets of combinatorial objects $\Omega$  have inner spaces of ``polynomials'' associated with them, called {\em polynomial spaces}. By ``polynomial'' is here meant a {\em zonal polynomial} with respect to a point $a$ in $\Omega$; this is a function on $\Omega$ obtained by composing a polynomial in one variable over some field $\mathbb F$ with the function $\rho(a,\cdot)$, where $\rho$ is a symmetric function from $\Omega\times \Omega$ to $\mathbb F$ (called a separation function). The most important example is that in which $\Omega$ is the set of vertices of a {\em distance-regular graph} (or of a generalization of it, called {\em association scheme}), and $\rho$ is essentially the distance function of the graph. The distance between two vertices of a graph  is the minimum number of edges in a path joining them; a connected graph is distance-regular if the number of vertices at any given distances $i$ from a vertex $u$ and $j$ from a vertex $v$ only depends on the distance between $u$ and $v$. A special case is the {\em Johnson graph}, whose vertices are the $k$-subsets of a finite set, and two $k$-subsets are joined by an edge if and only if their intersection has cardinality $k-1$. The most important application of  polynomial spaces is to deriving tight bounds for the sizes of combinatorial configurations with a high degree of symmetry in terms of the dimensions of certain subspaces of a polynomial space. The configurations corresponding to the Johnson graph are those called {\em designs}. A $t$-$(v,k,\lambda)$ design is a family of $k$-subsets (blocks) of a set $X$ of size $v$ with the property that all $t$-subsets of $X$  are contained in $\lambda$ blocks. The general technique mentioned above has applications to {\em coding theory} too.

{\bf 3. Combinatorics and geometry.} The main connection between combinatorics and geometry is related to {\em polytopes}. A $d$-polytope is a full-dimensional bounded intersection of closed subspaces in ${\mathbb R}^d$. For instance, the Platonic solids are 3-polytopes. Higher dimensional polytopes started to be investigated only in the 50's, due to their importance for {\em linear programming} (this refers to techniques for optimizing a linear functional subject to a collection of linear constraints). One defines the {\em faces} of a polytope as the areas of contact if the polytope is made to touch a plane surface in ${\mathbb R}^d$. There are only finitely many faces of each dimension $0,1,\ldots, d-1$, so one can form a vector by collecting these numbers; this is called the {\em f-vector} of the polytope. According to the generalized Euler formula, the alternating sum of the entries of an $f$-vector is either 0 or 2, depending on the parity of $d$. The important question arises what other relations, if any, do the entries of an $f$-vector satisfy. The best hope is to find a complete characterization of $f$-vectors. The answer is known in dimension 3 as a list containing Euler's formula plus two more inequalities. In general, there is a recent theorem which relates $f$-vectors of {\em simplicial polytopes} (i.e., all faces are simplices) to sequences whose entries count monomials of a certain degree in a certain collection of monomials in several variables. 

A related problem to the characterization of $f$-vectors of $d$-polytopes is the purely combinatorial characterization of the {\em  boundary complex} of a $d$-polytope (known as the {\em Steinitz problem}).  The boundary complex of a $d$-polytope is formed by its faces other than the polytope itself. 
Let us describe briefly what the Steinitz problem is about. To do this, we need to introduce {\em polyhedral complexes}, which  are defined as collections ${\cal C}$ of polytopes such that each face of a member of ${\cal C}$ is a member of ${\cal C}$, and the intersection of two members of ${\cal C}$ is a face of each; the dimension of ${\cal C}$ is the maximum of the dimensions of its members. A complex is called {\em polytopal} if it is isomorphic to the boundary complex of a polytope. The Steinitz problem asks for necessary and sufficient conditions for a complex to be polytopal. The structure of the  boundary complex of a polytope is encoded in the corresponding {\em face-lattice}, that is, the collection of all faces ordered by inclusion. Thus we also have the problem of finding geometric realizations of face-lattices. The Steinitz problem was completely solved in dimension 3, but not in higher dimensions. An important obstruction to a complex being polytopal is given by {\em shellability}, since it is known that each polytope is shellable. A pure complex of dimension $d$ (a complex is pure if the union of its members coincides with the union of its $d$-dimensional members) is defined to be shellable if there is an ordering of its members (faces), such that the intersection of each face with the union of the previous ones is a $(d-1)$-ball or a $(d-1)$-sphere. The difficulty of the Steinitz problem in dimensions higher than 3 shifted the focus of research to algorithmic approaches. For instance, it is known that there exists an algorithm for deciding whether a given complex is polytopal, and there are practical approaches in this direction. Some of the results mentioned above have been used in the enumeration of combinatorial types.

Another important research area is concerned with {\em point lattices}. By a point lattice is here meant a set of points in ${\mathbb R}^d$ which form a discrete additive subgroup. There are interesting combinatorial results concerned, for instance, with: 1) the relationship between the volume of a convex body on the one hand, and the number of lattice points it contains, or some invariants of the lattice; 2) covering and packing problems, that is, problems related to sets of translates of a convex body; 3) lattice polytopes, that is, polytopes all of whose vertices are in the integer lattice ${\mathbb Z}^d$.

{\bf 4. Combinatorial representation theory.} A representation of a group $G$ is a group homomorphism from $G$ to the group of invertible linear maps from a certain vector space to itself (i.e., the general linear group of the vector space). A representation is called {\em irreducible} if the corresponding vector space has no nontrivial invariant subspaces under the obvious action of $G$. The {\em character} of the representation is the map which associates with every element of $G$ the trace of the corresponding linear map. The main questions related to representations of groups are: 1) how do we index/count the irreducible representations? 2) what are their dimensions? 3) what are their characters? 4) how do we construct irreducible representations? 5) what are the multiplicities of irreducible representations in the direct sum decomposition of the tensor product of two irreducible representations? 6) how do we find interesting representations?

It turns out that, for the symmetric group $S_n$, the answer to most of these questions is given in terms of standard Young tableaux, which are semistandard Young tableaux with distinct entries (see section 1 for the definition of the latter). On the other hand, one can construct interesting representations for $S_n$ by considering actions on combinatorial structures, such as partially ordered sets (of partitions of a set, of labeled trees etc.), maximal chains in partially ordered sets, and the free Lie algebra.
In a similar way, the answers for the general linear group $GL_n({\mathbb C})$ are given in terms of semistandard tableaux of a given shape with entries $1,\ldots,n$. For instance, the irreducible polynomial representations are indexed by partitions with at most $n$ rows; the character of the representation indexed by $\lambda$ is a sum over semistandard Young tableaux $T$ of shape $\lambda$ of monomials obtained by associating with each entry $i$ of $T$ the variable $x_i$, and by taking the product of these variables. This character is a symmetric polynomial in $x_1,\ldots, x_n$, known as a {\em Schur polynomial}. The decomposition of tensor products of irreducible representations can thus be reduced to multiplying Schur polynomials and expressing the result as a linear combination of such polynomials. The corresponding nonnegative integer coefficients (known as {\em Littlewood-Richardson coefficients}) have a combinatorial interpretation in terms of  tableaux, which is known as the {\em Littlewood-Richardson rule}.

 The above results were generalized to {\em finite groups of Lie type} and to {\em modular representations}. They were also generalized to other {\em reflection groups} and {\em Lie groups}, as well as to their $q$-deformations known as {\em Hecke algebras} and {\em quantum groups}. The answers are no longer in terms of tableaux, but in terms of more general combinatorial structures, such as: 1) the chambers  of the complement in ${\mathbb R}^n$ of a certain collection of hyperplanes; 2) certain paths (piecewise-linear curves) in a given vector space, known as {\em Littelmann paths}; 3) lattice points in certain polytopes. Several other difficult combinatorial problems arise in this context. Some of them are related to {\em Kazhdan-Lusztig polynomials}; these are polynomials in $q$ with nonnegative integer coefficients which arise as entries of the transition matrix between two bases of the Hecke algebra of type $A$ (a $q$-deformation of the group algebra of the symmetric group).

The theory of symmetric functions, which is to a great extent combinatorial, was very much influenced by representation theory. A recent development in this area of interplay is the theory of {\em Macdonald polynomials}. An important breakthrough was the proof of the fact that a slight modification of these polynomials expand in the basis of Schur polynomials with coefficients which are polynomials in $q,t$  (the {\em Kostka-Macdonald polynomials}) with nonnegative integer coefficients. The proof uses sophisticated methods from algebraic geometry. Combinatorial interpretations for the coefficients of the Kostka-Macdonald polynomials exist only in special cases.

{\bf 5. Combinatorics and group theory.}  Like in other cases we mentioned, the relationship between group theory and combinatorics is a two-way relationship. Let us start by pointing out the way in which group theoretical tools were used in combinatorics. In fact, we have already mentioned one such instance, namely the enumeration of orbits of a finite group acting on a finite set. Indeed, the number of such orbits can be expressed in terms of the character of the corresponding permutation representation (see the definition of characters in section 4). Another application of group representations is related to the enumeration of solutions to equations in groups. To give a flavour of these results, consider the function associating with every element of the symmetric group the number of its factorizations with a fixed number of factors, in which each factor is a fixed power (possibly negative) of some permutation. It turns out that this function is a character of the symmetric group, which was determined in special cases.
Last but not least, we mention that group theoretical results, including the classification of finite simple groups, were used to derive combinatorial results, mostly concerned with {\em transitive actions} of finite groups (recall that a group acts transitively on a set if given any two elements of the set, there is an element of the group carrying the first to the second). For instance, such tools were used in the theory of {\em distance-transitive graphs}; these are graphs with the property that their automorphism group acts transitively on ordered pairs of vertices at any given distance. For instance, it was proved in this way that there are finitely many such graphs with vertices of any fixed degree (that is, a fixed number of neighbors). A project is underway to classify all distance-transitive graphs.

Combinatorial tools were also successfully used in group theory, including the classification of finite simple groups. The main point is that many interesting finite groups appear as automorphism groups of certain combinatorial structures with a considerable amount of symmetry. In fact, any finite group is the automorphism group of a finite {\em strongly regular graph}, and of a {\em Steiner triple system}. A graph is strongly regular if it is regular (all vertices have the same number of neighbors), not complete or empty, and the number of vertices adjacent to two vertices $s$ and $t$ only depends on whether $s$ and $t$ are adjacent or not; a Steiner triple system is a $2$-$(v,3,1)$ design (see section 2). Given the result above, it comes as no surprise that the most complex finite simple groups, namely the {\em sporadic} ones, are best described as automorphism groups of some of the previously mentioned structures. Furthermore, there are classification results related to the corresponding actions, for instance to the action of the sporadic groups on distance-transitive graphs. 

Another area of interplay between combinatorics and group theory is related to computing in permutation groups. Some computational problems, such as determining the center, the composition factors, or the Sylow subgroups of a group, are known to be polynomial (see the discussion about complexity in section 10). Others are NP-hard, and, in fact, are at least as hard in general as the problem of deciding if two graphs are isomorphic. Some of the problems in the latter category are: isomorphism of permutation groups, finding generators for the automorphism group of a graph, finding generators for the intersection of two permutation groups etc.

{\bf 6. Combinatorics and commutative algebra.} The main object of interest in this area is the {\em Stanley-Reisner ring} of a finite {\em simplicial complex} (i.e., a family of subsets of a set - called simplices - which is closed under taking subsets); this is defined as a quotient of a polynomial ring (with variables corresponding to the vertices of the complex) by the ideal generated by the squarefree  monomials corresponding to subsets of vertices which are not simplices. In particular, we have the Stanley-Reisner ring of a partially ordered set, which is constructed in terms of the simplicial complex of chains in the partially ordered set (known as the {\em order complex}). Clearly, the whole algebra structure of the Stanley-Reisner ring is determined by the combinatorial data in the simplicial complex or the partially ordered set. There are many ramifications from this observation. For instance, certain invariants of the Stanley-Reisner ring, such as its Hilbert function and certain vectors called {\em h-vectors}, are determined by the $f$-vectors of the simplicial complex (introduced in section 3 for polytopes); we note in passing that there are several results related to the characterization of both $h$- and $f$-vectors in certain cases. Furthermore, there are strong conditions on the combinatorics of a partially ordered set to allow for the {\em Cohen-Macaulay} property of the corresponding Stanley-Reisner ring. Recall that the Cohen-Macaulay property of a quotient of a polynomial ring by a homogeneous ideal can be reformulated as the existence of a certain direct sum decomposition of the quotient ring, known as  the {\em Hironaka decomposition}; essentially, this property expresses a strong regularity property of the associated variety. A partially ordered set is defined to be Cohen-Macaulay if the associated Stanley-Reisner ring is; there is also a purely combinatorial characterization. We point out that testing a partially ordered set $P$ for the Cohen-Macaulay property is very useful since it is known that it implies the same property for an important class of algebras, called {\em algebras with straightening law over P}. Furthermore, every Hironaka decomposition of the Stanley-Reisner ring of $P$ induces a Hironaka decomposition of the algebra since it is known that the two are isomorphic as vector spaces.

Let us now turn to conditions on a partially ordered set which imply the Cohen-Macaulay property. First of all, a necessary condition is that the partially ordered set is {\em ranked}, that is, all maximal chains have the same length (this enables one to define the rank of each element). Another condition is given in terms of the M\"{o}bius function, discussed in section 2. The M\"{o}bius function is called {\em alternating} if its sign (on an interval) is determined by the parity of the difference between the  ranks of the endpoints of the interval. It is known that the alternating property of the M\"{o}bius function implies the Cohen-Macaulay property. On the other hand, it is known that {\em lexicographic shellability} of the partially ordered set also implies the Cohen-Macaulay property, and one can derive the Hironaka decomposition of the Stanley-Reisner ring from a {\em shelling}. A partially ordered set is lexicographically shellable if it is ranked, and its set of maximal chains admits a linear ordering (called shelling) such that every chain contains a unique minimal subset that is not contained in any previous chain. Lexicographic shellability is a method for proving shellability (see section 3) of order complexes of partially ordered sets. An important class of partially ordered sets which is known to be (lexicographically) shellable is formed by all finite quotients of a Coxeter group by a parabolic subgroup, ordered by {\em Bruhat order}; the latter is defined on the Coxeter group by setting $u\le v$ if every {\em reduced decomposition} for $v$ (or word, that is, a factorization of $v$ of minimum length in terms of the Coxeter generators of the group) contains a subword which is a reduced decomposition for $u$. Similarly, the lattice of subgroups of a finite group is (lexicographically) shellable if and only if the group has a property known as {\em supersolvability}.

{\bf 7. Combinatorics and algebraic geometry.} There are several applications of results in algebraic geometry to combinatorics. For instance, known estimates for the number of connected components of real varieties were used to estimate the number of tests needed to solve the membership problem for certain subsets of ${\mathbb R}^d$ (i.e., to decide if a certain point is in the subset or not); a typical example is to decide if the numbers in a finite set are all distinct. Another example is using results about the number of solutions of polynomials with several variables over a finite field to prove that certain graphs contain regular subgraphs (see section 5). Algebraic relations among the subdeterminants of a matrix (the {\em Grassmann-Pl\"{u}cker relations}) were used to recognize polytopal complexes (see the discussion in section 3).  Finally, results in the theory of elliptic curves have been used to construct codes with excellent error-correcting capacity. 

As far as combinatorial methods in algebraic geometry are concerned, they are related to  combinatorial descriptions of certain algebraic varieties, including their cohomology and {\em singular locus} (recall that a point is called smooth if the dimension of the tangent space at that point is equal to the dimension of the variety; otherwise, it is called singular). 

A class of varieties where such methods were used successfully is that of {\em Schubert varieties}. They arise in a certain decomposition of varieties of flags in ${\mathbb C}^n$; the latter can be viewed as quotients of the group $SL_n$ of $n\times n$ complex matrices of determinant 1 by certain subgroups. If we consider the quotient $SL_n/B$ by the subgroup $B$ of upper triangular matrices, the Schubert varieties are indexed by permutations of $n$ elements. Essentially, the Schubert decomposition of $SL_n/B$ is the geometric counterpart of the well-known Gaussian elimination method in linear algebra. The {\em Grassmannian} of $k$-subspaces in ${\mathbb C}^n$ is a special case of a flag variety. It turns out that the number of points in a suitable triple intersection of Schubert varieties in Grassmannians is just a Littlewood-Richardson coefficient, discussed in section 4. Hence it is described by a combinatorial rule in terms of Young tableaux. However, there is as yet no good explanation for the fact that these numbers appear both in Schubert calculus and the representation theory of $GL_n({\mathbb C})$. On the other hand, the Littlewood-Richardson coefficients have found a surprising application recently to describing eigenvalues of sums of Hermitian matrices. The problem of finding a combinatorial interpretation for the intersection numbers of Schubert varieties in arbitrary  flag varieties (which is equivalent to describing the multiplicative structure of the ordinary cohomology of flag varieties with respect to a certain natural basis) is not yet solved; however, it seems that the answer should be in terms of certain diagrams of strands (known as {\em rc-graphs}), generalizing Young tableaux. Another problem for Schubert varieties with combinatorial ramifications which has been investigated is the membership problem.

We now turn to the problem of the singular locus of Schubert varieties in $SL_n/B$. It is known that the Schubert variety indexed by a permutation $w$ is singular if and only if $w$ contains elements which relate to each other in the same way as the elements of the permutations (patterns) 4231 and 3412 do. Recently, a purely combinatorial description of the irreducible components of the singular locus of the previously mentioned Schubert varieties was given. On the other hand, there is a nice connection between singular points in Schubert varieties and the Kazhdan-Lusztig polynomials discussed in section 4; essentially, a point is smooth if and only if the corresponding Kazhdan-Lusztig polynomial is equal to 1. This connection relates the combinatorics in the two areas.

The previous results were generalized to some extent to Schubert varieties corresponding to other Lie groups than $SL_n$, as well as to other varieties, such as {\em toric} and {\em quiver varieties}. On the other hand, more complex cohomology theories have been studied, such as {\em equivariant}, {\em quantum}, and {\em intersection cohomologies}, {\em K-theory} etc. Let us just mention that the analogues in quantum cohomology of the intersection numbers of Schubert varieties discussed above are the so-called {\em Gromov-Witten invariants}, enumerating certain curves in $SL_n/B$. There are some recent combinatorial results related to them.

{\bf 8. Combinatorics and topology.} The following idea was successfully used in combinatorics: associate a topological space with a combinatorial object, and derive combinatorial results from topological ones applied to this space; furthermore, the topological study of the space is interesting in its own right. One can associate topological spaces with partially ordered sets, graphs, and matroids (the latter are abstract generalizations of matrices). For instance, by applying topological results to the simplicial complex of chains in a partially ordered set (its order complex), such as fixed point theorems or homology calculations, one can obtain fixed point theorems and other results for the partially order set (such as results about the M\"{o}bius function, introduced in section 2, which is just the reduced {\em Euler characteristic} in this context). There are related results referring to the homotopy type of a partially ordered set (that is, of its order complex). For instance, it is known that every lexicographically shellable partially ordered set (such as the Bruhat order on Coxeter groups, see section 6)  has the homotopy type of a wedge of spheres. Shellings are also useful for constructing bases for the cohomology of the order complex. In some cases (such as partially ordered sets of partitions of a set of various types), the symmetric group acts on the partially ordered set and hence on its cohomology, and one can use the previously mentioned constructions to identify the representation. Finally, we mention results related to concrete realizations of the order complex as a finite {\em CW complex} (also known as a cell complex). Such a construction, where in addition the boundaries of cells have the homology of spheres, was recently given for intervals in the Bruhat order on the symmetric group.

Another area of combinatorics related to topology is {\em topological graph theory}, which studies embeddings of graphs in surfaces. It turns out that the graphs which cannot be embedded in a given 2-dimensional surface can be described in terms of a finite list of nonembeddable graphs, called {\em forbidden minors}. For instance, there are only 2 forbidden minors for the plane, but 35 of them for the projective plane.

Conversely, combinatorial methods are useful in topology. One such instance is related to knot theory. Indeed, one can construct knot invariants, such as the Jones polynomial, purely combinatorially,   by associating a graph with a knot, and then by considering polynomial invariants of the graph, such as the {\em Tutte polynomial}. Another instance of combinatorial methods in topology is related to the representation of certain manifolds by edge-colored graphs (known as {\em crystallizations} or, more generally, {\em combinatorial maps}); this representation led to several classification results.

Another important interaction between combinatorics and topology is related to {\em complex subspace arrangements}. These are finite sets of subspaces in a complex vector space. Special cases are hyperplane arrangements, for which all subspaces have codimension 1. An important combinatorial structure associated with the arrangement is the {\em intersection lattice}, that is, the set of intersections of subspaces partially ordered by reverse inclusion. Much research has been devoted to understanding the relationship between the intersection lattice and other combinatorial data on the one hand, and the topology of the complement of the arrangement, on the other hand (for instance, its cohomology, or the homotopy type of its intersection with the unit sphere, which is known to be that of a wedge of spheres in the case of hyperplane arrangements).

An even bolder step than the ones described above was taken by developing combinatorial versions of certain topological theories, such as {\em Morse theory} (this theory describes the relationship between a function's critical points and the homotopy type of the function's domain). One application of discrete Morse theory was to shellability  of simplicial complexes (see section 3). Another interesting recent theory is that of {\em combinatorial differential manifolds}; the key ingredient here is a discrete version of the tangent bundle, defined in terms of a combinatorial structure known as {\em oriented matroid}.

{\bf 9. Combinatorics and probability theory.} One instance of the interplay between combinatorics and probability theory is the use of probabilistic methods to prove the existence of certain combinatorial configurations. The idea is very simple: one creates a probability space and shows the existence of a good point (configuration) by showing that the bad points have measure less than one. Essentially, this method is based on the use of a relatively simple probabilistic principle, combined with the analysis of the asymptotic behavior of some combinatorial function. This method was applied to problems such as: 1) finding lower bounds for {\em Ramsey numbers} (a Ramsey number is essentially the minimum number $n$ such that all configurations of a given type with size at least $n$ have a given property); 2) problems in {\em extremal set theory} (this is concerned with finding the maximum size of a family of subsets of a set having a certain property); 
3) the probabilistic analysis of algorithms (see section 10).

An important area of interaction between probability theory and combinatorics is that of {\em random graphs}. A random graph is a family of graphs endowed  with a probability distribution. A very popular example is the {\em binomial random graph} $G(n,p)$, in which there are $n$ vertices and each edge is chosen with probability $p$. An important problem is to study the evolution of $G(n,p)$ as $p$ increases from 0 to 1. It turns out that, up to a certain known threshold, $G(n,p)$ consists almost surely of components with at most one cycle, after which it almost surely consists of a unique giant component and small components with at most one cycle. This phenomenon is called {\em phase transition}. There are several estimations of the size of the giant component in the neighborhood of phase transition. There is a similar phenomenon to phase transition for certain graph properties, in the sense that the probability of $G(n,p)$ having the property jumps from 0 to 1 in the limit. Another branch of the theory of random graphs is concerned with limit distributions of numerical characteristics of a graph, such as the number of copies of a fixed graph it contains and its chromatic number (see the definition in section 2). 

Another area of interaction between combinatorics and probability theory we are mentioning is related to {\em random walks} on graphs. These are processes which start at a certain vertex of the graph, and at each step move to a neighboring vertex with a given probability. Examples include shuffling a deck of cards and random walks on groups, in which each move means multiplication by certain elements in the group. Assuming that the graph is connected and $d$-regular (all vertices have $d$ neighbors), that all transitions occur with probability $1/d$, and that we follow the random walk long enough, the current point will be almost surely uniformly distributed over the vertices. The analysis of random walks includes computation of the {\em mixing rate}, which measures how fast this convergence happens. This depends on global connectivity properties of the graph, such as expansion (see the discussion in section 2). In the case of random walks on groups, the main tools for estimating the mixing rate come from representation theory. An important application of random walks is to {\em randomization of algorithms}. Indeed, a useful heuristics for computationally hard problems involving local search is to randomly select a neighbor of the current state at each step. This allows the process to reach ``obscure'' parts of large sets and avoid ``tricky'' corners. Concrete applications include estimation of the volume of a bounded convex body, enumeration of matchings in graphs, statistical sampling etc.

Finally, we would like to mention a relatively recent development, namely {\em free probability} theory. This is an analog of classical probability in which the random variables are replaced by elements of some (possibly noncommutative) algebra. The central notion is an analog of the notion of independence for random variables. There are some subtle connections between this theory and combinatorics, particularly combinatorics related to the representation theory of the symmetric group (see section 4). 

{\bf 10. Combinatorics and theoretical computer science.} Many algorithms in computer science are concerned with combinatorial structures, such as graphs. Assuming that the edges are weighted, here is a list of computational problems involving graphs: 1) find the tree of minimum total weight containing all vertices (the {\em minimum spanning tree}); 2) viewing the graph as a network of pipes and the edge weights as capacities, find the maximum flow in the network; 3) find the path containing all vertices and minimizing the sum of weights on its edges (the {\em traveling salesman problem} TSP); 4) find a maximum cut in the graph, that is a subset of the vertices maximizing the sum of weights of edges with exactly one endpoint in the subset (the {\em max-cut problem}). A major theme in computer science is estimating the {\em complexity} of various algorithms for solving such problems. By complexity is here meant the amount of time (or steps) and computational resources needed. By constructing algorithms, one shows that a task can be done in a certain number of steps. It is often the most difficult part to show that there is no faster way, i.e. requiring fewer steps. Surprisingly, there is a nice topological approach to this problem in some cases; essentially, one associates a topological space with the given problem, and observes that the the more complicated this space is topologically, the higher the complexity of the algorithm.  For a large class of problems, there are {\em polynomial algorithms} for solving them (that is, their complexity is a polynomial in the number of vertices of the graph, say). For other problems, called {\em NP-hard}, no such algorithm is known. For instance, problems 1) and 2) are polynomial, while 3) and 4) are NP-hard. Membership in the latter class is shown by reducing the given problem to a standard one, such as TSP. If this is the case, one looks for heuristics for finding an approximate solution. For example, the best known heuristics for TSP has a performance ratio of $3/2$; this means that the weight of the path found by the heuristics is at most $3/2$ as large in general as the minimum weight.

Here are some examples of heuristics for approximating solutions to NP-hard problems: 1) the greedy heuristics, based on making the best choice at each step of the algorithm; 2) local improvement techniques, such as randomization (see the discussion in section 9); 3) heuristics based on linear programming (see the discussion in section 3 about this method). The mathematical formalization of some of these heuristics, such as the greedy one, led to combinatorial structures with many interesting properties, such as those called {\em greedoids}. 

The difficulty of NP-hard problems led to the development of new models of computation. We enumerate some of them: 1) {\em parallel computation}, that is, algorithms in which several tasks are performed simultaneously; 2) {\em on-line algorithms}, that is, algorithms in which not all data is available at once, but becomes available as the algorithm progresses; 3) {\em probabilistic analysis of algorithms}, i.e. an average-case analysis based on random input (as opposed to the classical worst-case analysis). For instance, there is a heuristics for TSP which, under probabilistic analysis, has a performance ration of $1/\sqrt{n}$, where $n$ is the number of vertices of the graph. The mathematical formalization of the above models of computation is only in its initial stages.

\vspace{1cm}

Our list of areas of mathematics and computer science interacting with combinatorics is by no means exhaustive. For instance, we have not mentioned {\em combinatorial number theory}, concerned with the study of properties of ``dense'' subsets of positive integers, representations of positive integers as sums of elements in a given set, congruences for combinatorial numbers (such as the binomial coefficients) etc. On the other hand, we have not devoted a special section to areas such as analysis and differential equations (the connections between the latter and combinatorics are just starting to emerge); however, we have mentioned structures and tools from analysis used in combinatorics (for instance, generating functions and polynomial spaces in sections 1 and 2, as well as the asymptotic behavior of combinatorial functions, in section 9). Finally, we have not mentioned applications of combinatorics to other areas of science. Such areas include: operations research, electrical engineering, statistical physics, chemistry, and molecular biology.

\end{document}